\documentclass[letterpaper,10pt]{article}
\pdfoutput=1
\usepackage{mathtools}
\usepackage{amssymb}
\usepackage{mathrsfs,dsfont,textcomp}

\usepackage{nicematrix}

\usepackage{empheq}

\usepackage{amsthm}
\usepackage{thmtools} 

\usepackage{titlesec} 

\usepackage{parskip}
		
\usepackage{enumitem}

\usepackage{eso-pic}
\usepackage{tikz}
\usetikzlibrary{tikzmark}

\usepackage{footnote}
\usepackage{footmisc} 
\newcounter{savefootnote}
\newcounter{symfootnote}
\newcommand{\symfootnote}[1]{%
   \setcounter{savefootnote}{\value{footnote}}%
   \setcounter{footnote}{\value{symfootnote}}%
   \ifnum\value{footnote}>8\setcounter{footnote}{0}\fi%
   \let\oldthefootnote=\thefootnote%
   \renewcommand{\thefootnote}{\fnsymbol{footnote}}%
   \footnote{#1}%
   \let\thefootnote=\oldthefootnote%
   \setcounter{symfootnote}{\value{footnote}}%
   \setcounter{footnote}{\value{savefootnote}}%
}

\usepackage{etoolbox}
\patchcmd{\thebibliography}{\section*{\refname}}{}{}{}

\usepackage[letterpaper, top=3cm, bottom=2cm, left=2.5cm, right=2.5cm]{geometry}

\titleformat
{\section}
{\scshape} 
{\thesection.} 
{1.0pc} 
{} 
\titlespacing  
{\section}
{0pc} 
{1.5pc} 
{\bigskipamount} 

\titleformat
{\subsection}
{\bfseries}
{\thesubsection}
{1.0pc}
{}
\titlespacing
{\subsection}
{0pc}
{\bigskipamount}
{\medskipamount}

\declaretheoremstyle[
spaceabove=\medskipamount, spacebelow=\medskipamount,
headfont=\normalfont\bfseries,
bodyfont=\normalfont,
postheadspace=  12pt,]
{mystyledefinition}
\declaretheorem[style=mystyledefinition, numberwithin=section]{definition}

\declaretheoremstyle[
spaceabove=\medskipamount, spacebelow=\medskipamount,
headfont=\normalfont\sl,
bodyfont=\normalfont,
postheadspace=  12pt,]
{mystyleexample}

\declaretheoremstyle[
spaceabove=\medskipamount, spacebelow=\medskipamount,
headfont=\normalfont\bfseries,
bodyfont=\sl,
postheadspace=  12pt,  
]{mystylethm}
\declaretheorem[style=mystylethm, numberwithin=section, numberlike=definition]{theorem}

\declaretheorem[style=mystylethm, numberwithin=section, numberlike=definition]{lemma}

\font\titlefont= cmbx12
\font\authorfont=cmcsc10 
\font\headfont=cmcsc10

\font\small=cmr9
\font\smallbf=cmbx9

\font\smalltt=cmtt9

\def\Proof{\noindent {\sl Proof.}\enspace}

\def\qedmark{\hbox{\vrule height 4pt width 3pt}}
\def\qedskip{\vrule height 4pt width 0pt depth 1pc}
\def\qed{\penalty 1000\quad\penalty 1000{\qedmark\qedskip}}

\def\cA{{\cal A}}

\def\cF{{\cal F}}

\def\cO{{\cal O}}

\def\cU{{\cal U}}
\def\cV{{\cal V}}

\def\srC{\mathscr{C}}

\def\gF{\mathfrak{F}}

\def\ovQ{\overline{Q}}

\renewcommand{\r}{\rightarrow}

\begin{document}

\begin{center}
{\titlefont The Spectral Flow of a Restriction  to a Subspace and \par the Maslov Indices in a Symplectic Reduction \footnote{{\bf MSC(2020).} 47A53, 53D12, 58B15.}}
\end{center}
\bigskip
\centerline{{\authorfont Henrique Vit\'orio}\symfootnote{Departamento de Matem\'atica, 
Universidade Federal de Pernambuco, Brazil. E-mail address: \smalltt henrique.vitori@ufpe.br} }

\vskip .5cm

\bigskip
\medskip
\hfil {\hsize = 10cm \vbox{\noindent {\smallbf Abstract.} \small
We give a simple proof of a known formula that relates the spectral flow of a continuous path of quadratic forms
of Fredholm type with the spectral flow of the restrictions of the forms to a fixed closed finite 
codimensional subspace. We then apply this to obtain a formula relating the Maslov index of a continuous path
in a Fredholm Lagrangian Grassmannian with the Maslov index of its symplectic reduction by a closed
finite codimensional coisotropic subspace.}} \hfil

\vskip.3cm

\bigskip

\noindent 

\section{Introduction}

\noindent Given a real finite dimensional vector space $E$, a quadratic form $Q$ on $E$ and a vector subspace $W\subset E$, the index 
of $Q$ differs from the index of its restriction to $W$ by the following expression
\begin{equation}\label{eq1}
{\rm ind}\, Q - {\rm ind}\, Q|_W = {\rm ind}\, Q|_{W^{\perp_Q}} + \dim(W\cap W^{\perp_Q}) - \dim(W\cap\ker Q),
\end{equation}
where $W^{\perp_Q}=\left\{ X\in E \,\; | \,\; Q[X,Y]=0 \;\; \forall Y\in W \right\}$. 
This linear algebra fact played a key role in the establishing in \cite{ballmann}
of a Morse index theorem for closed Riemannian geodesics from the theorem for geodesic segments. 
When dealing with strongly indefinite problems (e.g. when working on a semi-Riemannian background), the notion of Morse index 
is replaced by the notion of {\sl spectral flow} (see Sec. \ref{sectionspectralflow}), and it becomes relevant to extend formula (\ref{eq1}) to a 
formula relating spectral flows. In \cite{benevieri} the authors proved the following:

\begin{theorem}\label{maintheorem}
Let $Q=\{Q_t\}_{t\in [a,b]}$ be a continuous path of quadratic forms of Fredholm type on a separable infinite dimensional
real Hilbert space $H$. Let $V\subset H$ be a closed finite codimensional subspace and denote by $\overline{Q}=\{\overline{Q}_t\}$
the path of quadratic forms on $V$ obtained by restricting the $Q_t$ to $V$. Then the spectral flows ${\rm sf}(Q,[a,b])$ and 
${\rm sf}(\overline{Q},[a,b])$ are related by
\begin{equation}\label{eq2}
\begin{aligned}  
{\rm sf}(Q,[a,b]) - {\rm sf}(\overline{Q},[a,b]) &  \; = \;  {\rm ind}\, Q_a|_{V^{\perp_{Q_a}}} - {\rm ind}\, Q_b|_{V^{\perp_{Q_b}}} 
+ {\rm dim} (V\cap V^{\perp_{Q_a}}) - {\rm dim} (V\cap V^{\perp_{Q_b}}) \\
& \;  \quad - \dim (V \cap \ker Q_a) + \dim (V\cap \ker Q_b).
\end{aligned}
\end{equation}
\end{theorem}

It turns out that being able to relate the spectral flows of $\{Q_t\}$ and $\{\overline{Q}_t\}$ plays a key role in many places. For instance,
\begin{itemize}
\item[1--] In the work \cite{benevieri} the authors used formula (\ref{eq2}) to
work out a Morse index theorem for closed semi-Riemannian geodesics from the same theorem for
semi-Riemannian geodesic segments;     

\item[2--] For establishing in \cite{vitoriobifurcation} a Morse index theorem for electromagnetic geodesic segments constrained to a fixed energy level, 
one had to resort to such a formula;

\item[3--] When performing a symplectic reduction of a path
$\ell:[a,b]\r \cF\Lambda_{L_0}$ in a Fredholm Lagrangian Grassmannian (see Sec. \ref{sectionmaslov}) by a closed finite codimensional 
coisotropic subspace $\mathds{W}\subset\mathds{H}$ of a symplectic Hilbertable space, formula (\ref{eq2}) leads naturally to a formula relating the Maslov index $\mu_{L_0}(\ell,[a,b])$ of $\ell$ with 
the Maslov index of the reduction of $\ell$. Indeed, one of the purposes of these notes is to prove the following: 
\end{itemize} 

\begin{theorem}\label{thmreduction}
Let $L_0\in \Lambda(\mathds{H})$ have null intersection with $\mathds{W}^{\perp_\omega}$ and let
$\ell:[a,b]\r\cF\Lambda_{L_0}(\mathds{H})$ be a continuous path with $\ell(t)\cap\mathds{W}^{\perp_\omega}=\{0\}$ for all $t$.
Set $\overline{L}_0=q(L_0\cap\mathds{W})$, where $q:\mathds{W}\r \mathds{W}/\mathds{W}^{\perp_\omega}$ is the quotient map. Then
$t\mapsto \overline{\ell}(t)=q(\ell(t)\cap\mathds{W})$ defines a continuous path in 
$\cF\Lambda_{\overline{L}_0}(\mathds{W}/\mathds{W}^{\perp_\omega})$ and
\begin{equation}\label{eq3}
\begin{aligned}
\mu_{L_0}(\ell, [a,b]) - \mu_{\overline{L}_0}(\overline{\ell},[a,b]) & \; = \; 
{\rm ind}\, \mathfrak{q}_a - {\rm ind}\, \mathfrak{q}_b + \dim \pi\bigl( \ell(a) \cap (V + \mathds{W}^{\perp_\omega}) \bigr)
- \dim \pi\bigl( \ell(b) \cap (V + \mathds{W}^{\perp_\omega}) \bigr) \\
& \; \quad - \dim (\ell(a)\cap V) + \dim (\ell(b)\cap V).
\end{aligned}
\end{equation}
In the above,
\begin{enumerate}[left=4pt .. \parindent]
\item $V=L_0\cap\mathds{W}$, so $V$ is a closed finite codimensional subspace of $L_0$;
\item $\pi: L_0+\mathds{W}^{\perp_\omega} \r L_0$ is the projection;
\item for $t=a,b$, $\mathfrak{q}_t$ is the quadratic form on the finite dimensional space $E_t=\pi\bigl(\ell(t)\cap(L_0+\mathds{W}^{\perp_\omega})\bigr)$
given by 
\begin{equation}\nonumber
\mathfrak{q}_t[u,v] = \omega\Bigl[ \left(\pi\big|_{E_t}\right)^{-1} u , v \Bigr].
\end{equation}
\end{enumerate}
\end{theorem}

The ideia of applying a symplectic reduction when dealing with the Maslov index has been applied in many 
contexts. For instance, it has been explored in the work \cite{maslovbanach} to obtain an intrinsic reduction to finite dimensions 
when working out a construction of the Maslov index in symplectic Banach spaces (see also \cite[Prop. 1.15]{Nicolaescu} 
and \cite{furutani2} for finite dimensional reductions in the context of symplectic Hilbert spaces). Let us also
recall that Booss-Bavnbek and Furutani (\cite{furutani2}) gave a definition of an infinite dimensional {\sl Hörmander index} 
$\sigma_{\small\mbox{Hör}}(L_1,L_2;L_0',L_0)$ for a quadruplet of Lagrangians satisfying $\dim (L_0/L_0\cap L_0')<\infty$ and
$(L_0,L_1)$ and $(L_0,L_2)$ are Fredholm pairs. It turns out that the expression in the right-hand side of (\ref{eq3}) 
gives the Hörmander index for the quadruplet $(\ell(a),\ell(b); L_0', L_0)$, where $L_0'=L_0\cap\mathds{W} + \mathds{W}^{\perp_\omega}$.
A proof of this fact can be obtained from Theorem \ref{thmreduction} and from an equality  
$\mu_{L_0'}(\ell,[a,b])=\mu_{\overline{L}_0}(\overline{\ell},[a,b])$ (for $L_0'$ as before) that one gets by using similar
arguments as the ones used in the proof of (\ref{eq3}).  

Given the utility of formula (\ref{eq2}), it seems useful to have a simple proof of it. Let us recall that the original proof in
\cite{benevieri} started by first assuming that the path $\{L_t\}$ of self-adjoint
operators representing $\{Q_t\}$ is a perturbation of a {\sl fixed} symmetry $\mathfrak{I}\in GL(H)$ by a path of compact 
operators, in which case the spectral flow of $\{Q_t\}$ depends only on the endpoints $Q_a$ and $Q_b$ and is given by a difference
of relative dimensions of commensurable pairs of subspaces (see \cite[Thm. 4.4]{benevieri});
the general case then follows from the existence of cogredient parametrices (\cite[Thm. 2.1]{recht}) and a further orthogonal 
trivialization argument (\cite[Prop. 4.14]{benevieri}). It turns out that a much simpler argument can be used to prove 
formula (\ref{eq2}) and that shows in a clear way how the latter reduces to the finite dimensional formula (\ref{eq1}) by exploiting 
the invariance properties of the spectral flow. That argument, which avoids the technicality of relative dimensions and parametrices,
was employed in \cite{vitoriobifurcation} by assuming non-degeneracy of the paths 
$\{Q_t\}$ and $\{\overline{Q}_t\}$ at the endpoints, and we shall extend it here so as to include degenerate endpoints.

This work is organized as follows. In Sec. \ref{sectionpreliminaries} we recall some concepts and prove some results about bilinear forms (symmetric and skew-symmetric)
of Fredholm type. After collecting some properties of the spectral flow in Sec. \ref{sectionspectralflow}, we prove Theorem \ref{maintheorem} 
in Sec. \ref{sectionproof}. In Sec. \ref{sectionmaslov} we begin by recalling the concepts surrounding the Fredholm Lagrangian Grassmannian
and the Maslov index of paths. The proof of Theorem \ref{thmreduction} is then carried out in Sec. \ref{sectionreduction} 
and Sec. \ref{sectionproofreduction}.

\section{Some preliminaries}\label{sectionpreliminaries}

\noindent Throughout this section, $H$ shall denote a separable infinite dimensional real {\sl Hilbertable} space. Let 
$Q$ be a bounded bilinear form on $H$ which is either symmetric (in which case we shall call it a quadratic form) or 
skew-symmetric. For a subspace $W\subseteq H$, we set $W^{\perp_Q}=\{ u \in H \,\; | \,\; Q[u,w]=0 \;\; \forall w\in W\}$. 
Recall that $Q$ is said non-degenerate if $\ker Q = H\cap H^{\perp_Q}=\{0\}$, and that a subspace $W\subset H$ is
non-degenerate with respect to $Q$ if $\ker Q|_W = W\cap W^{\perp_Q}=\{0\}$. 
The form $Q$ is said of Fredholm type if the bounded self-adjoint (or skew-adjoint) operator that 
represents it with respect to some (and hence to all) compatible Hilbert inner product on $H$ is Fredholm. It is clear 
that if $Q$ is of Fredholm type and if $W\subset H$ is a closed finite codimensional subspace, then the restriction $Q|_W$ is
of Fredholm type and $W^{\perp_Q}$ is finite dimensional. Also, if $Q$ is non-degenerate and of Fredholm type, then $Q$ 
is strongly non-degenerate, i.e. the musical map $Q^\sharp:H \r H^*$ is an isomorphism.

For the next two lemmas, let $Q$ be a fixed quadratic/skew-symmetric form of Fredholm type on $H$.

\begin{lemma}\label{lemmaperpperp}
Let $U\subset H$ be a subspace that is either finite dimensional or closed and finite codimensional. Then,
$(U^{\perp_Q})^{\perp_Q} = U + \ker Q$.
\end{lemma}
\Proof
Let $\pi:H\r H/\ker Q$ be the quotient map and $\ovQ$ be the 
form on $H/\ker Q$ induced by $Q$. Then,
\begin{enumerate}
\item $\bigl(W^{\perp_{\ovQ}}\bigr)^{\perp_{\ovQ}}=W$
for any closed subspace $W\subset H/\ker Q$ since the musical map $\ovQ^\sharp:H/\ker Q \r (H/\ker Q)^*$ is an isomorphism;
\item $\pi(V^{\perp_Q})=\pi(V)^{\perp_{\ovQ}}$ for any subspace $V\subset H$, by a purely algebraic argument;  
\item $\pi(U)$ is closed since the restriction $\pi|_U:U\r H/\ker Q$ is a Fredholm
operator if $U$ is closed and ${\rm codim}\,U<\infty$ (since $\pi$ is already Fredholm) and Fredholm operators have closed range.  
\end{enumerate}
Therefore we can apply 1. to $W=\pi(U)$ and iterate 2. to obtain $\pi\bigl((U^{\perp_Q})^{\perp_Q}\bigr) = \pi(U)$. From
this equality and the fact that $\ker Q\subset(U^{\perp_Q})^{\perp_Q}$, we conclude the result.  
\qed

\begin{lemma}\label{lemmaalgebraic}
Let $W\subset V$ be closed subspaces of $H$ with finite codimension and with $W$ non-degenerate. Then,
$(W^{\perp_Q}\cap V)^{\perp_Q} = W\oplus V^{\perp_Q}$.
\end{lemma}
\Proof
Fist of all, observe that $E=W\oplus V^{\perp_Q}$ is closed in $H$: indeed, if $\pi:H\r H/W$ is the quotient map, then
$\pi$ is continuous, $\pi(E)$ is finite dimensional (hence, closed) and $W=\pi^{-1}(\pi(E))$.
Now, it is clear that $(W\oplus V^{\perp_Q})^{\perp_Q}=W^{\perp_Q}\cap (V^{\perp_Q})^{\perp_Q}$. Thus, since, 
by the previous lemma, $(V^{\perp_Q})^{\perp_Q}=V+\ker Q$, and since $\ker Q\subset W^{\perp_Q}$, we obtain
$(W\oplus V^{\perp_Q})^{\perp_Q}=W^{\perp_Q}\cap V + \ker Q$. On the other hand, since $W^{\perp_Q}\cap V$ is 
finite dimensional, Lemma \ref{lemmaperpperp}
also applies to show that $\left( (W^{\perp_Q}\cap V)^{\perp_Q} \right)^{\perp_Q} = W^{\perp_Q}\cap V + \ker Q$.
Therefore, $\left( (W^{\perp_Q}\cap V)^{\perp_Q} \right)^{\perp_Q}=(W\oplus V^{\perp_Q})^{\perp_Q}$. Computing the 
$\perp_Q$ of both sides of last equality via Lemma \ref{lemmaperpperp}, we obtain 
$(W^{\perp_Q}\cap V)^{\perp_Q} + \ker Q = W\oplus V^{\perp_Q} + \ker Q$. The result now follows since 
$\ker Q$ is contained in both $W^{\perp_Q}\cap V$ and $W\oplus V^{\perp_Q}$.
\qed

For the proof of Theorem \ref{maintheorem} we shall need to guarantee the existence of a closed finite codimensional
subspace that is non-degenerate with respect to {\sl two} quadratic forms (see Lemma \ref{lemmaW} of Sec. \ref{sectionproof}). 
The easiest way to do so seems to be establishing a much stronger result about the distribution of non-degenerate subspaces: 

\begin{lemma}\label{lemmadensity}
Let be fixed a quadratic form $Q$ of Fredholm type on $H$. The non-degenerate closed subspaces $W\subset H$ of a fixed finite codimension $k$ form an open subset of the 
Grassmann manifold $G^k(H)$ of all $k-$codimensional closed subspaces of $H$. If $Q$ is non-degenerate, then that subset is
also dense.
\end{lemma}
\Proof
We shall also denote by $G_k(H)$ the Grassmannian of all $k-$dimensional subspaces of $H$. Let us first prove the assertion about density.
This is equivalent to proving the density in $G_k(H)$ of the 
subset $\Sigma_k$ of non-degenerate subspaces. For, the non-degeneracy of $Q$ implies that a subspace $W\in G^k(H)$ is non-degenerate if, and only if, 
$W^{\perp_Q}$ is non-degenerate (this follows from Lemma \ref{lemmaperpperp}) and that the map $G^k(H)\r G_k(H)$, $W\mapsto W^{\perp_Q}$,  
is a homeomorphism. The proof of the density of $\Sigma_k$ in $G_k(H)$ is an application of the infinite dimensional Transversality Density Theorem (TDT)
(\cite[Thm. 19.1]{transversalmappings}) as we now explain. Let $\srC=\{X\in H\backslash\{0\} \;\, | \;\, Q[X]=0\}$ be the isotropic cone of $Q$.
Then a subspace $W\subset H$ is degenerate if, and only if, $W$ is somewhere tangent to $\srC$, i.e. there exists $X\in W\cap \srC$ such that
$W$ is contained in the tangent space $T_X\srC$ (this follows since $T_X\srC=\langle X\rangle^{\perp_Q}$).
Passing to the projective space $\mathds{P}(H)$, which is a smooth Hilbert manifold, that condition is equivalent to 
the (finite dimensional) projective subspace $\mathds{P}(W)\subset \mathds{P}(H)$ being somewhere tangent to the projective ``hiperquadric''
$\mathds{P}(\srC) \subset \mathds{P}(H)$, which is a one codimensional smooth submanifold of $\mathds{P}(H)$. 
In order to apply TDT, we need to put the set of embeddings 
$\left\{\mathds{P}(W)\hookrightarrow \mathds{P}(H) \;\, | \,\; W\in G_k(H) \right\}$ in a parametric form (at least in a neighborhood of each point of $G_k(H)$).
For this, given $W_0\in G_k(H)$, let $L(W_0,W_0^\perp)$ be the Hilbertable space of all (necessarily bounded) 
linear maps $W_0 \r W_0^\perp$. The set of $W\in G_k(H)$ with $W\cap W_0^\perp=\{0\}$ forms a neighborhood of $W_0$ on which a 
smooth parametrization $L(W_0,W_0^\perp)\r \cU_{W_0}$ is given by $f\mapsto {\rm graph}\, f=\{X+f(X) \,\; | \,\; X\in W_0\}$. We then define 
\begin{equation}\nonumber
F : L(W_0,W_0^\perp) \times \mathds{P}(W_0) \r \mathds{P}(H), \quad F(f,[X]) = [X+f(X)],
\end{equation}   
where $[-]$ stands for the projective class of a vector. The map $F$ is transversal to $\mathds{P}(\srC)$ since it is a submersion, hence TDT
applies to guarantee that
the set of $f$ for which $\mathds{P}({\rm graph}\, f)\hookrightarrow \mathds{P}(H)$
is transversal to $\mathds{P}(\srC)$ is dense. Therefore, the set of $W\in G_k(H)$ for which $\mathds{P}(W)\hookrightarrow\mathds{P}(H)$ 
is transverse to $\mathds{P}(\srC)$ is dense. 
\par Let us now consider the assertion about openness. Given $W_0\in G^k(H)$ with $W_0\cap W_0^{\perp_Q}=\{0\}$, there exist neighborhoods 
$\cA$ of $W_0$ in $G_k(H)$ and $\cA'$ of $W_0^{\perp_Q}$ in $G^k(H)$ such that $W\cap W' =\{0\}$ for all $W\in\cA$ and $W'\in\cA'$. 
On the other hand, since $Q$ may be degenerate, the correspondence $W \mapsto W^{\perp_Q}$ is not continuous on the whole of $G_k(H)$ (indeed, $\dim W^{\perp_Q}$ might jump). Nevertheless, one can show without difficulty that on the neighborhood $\cV$ of $W_0$ given by
$\cV=\{W \,\; | \,\; W\cap \ker Q= \{0\}\}$ the map $W\mapsto W^{\perp_Q}$ takes values in $G_k(H)$ and is 
continuous. Thus, for $W$ sufficiently close to $W_0$ we will have $W^{\perp_Q}\in\cA'$. This concludes the proof.     
\qed

\section{The spectral flow}\label{sectionspectralflow}

\noindent Let $H$ be a separable infinite dimensional real Hilbert space, and let
 $\hat{\gF}=\hat{\gF}(H)$ be the space of self-adjoint Fredholm operators 
on $H$. The space $\hat{\mathfrak{F}}$ has three connected components: two contractible components $\hat{\gF}_+$ and $\hat{\gF}_-$, 
which consist on the operators in $\hat{\gF}$ that are, respectively, essentially positive and essentially negative 
(i.e. that have finite dimensional negative, respectively, positive, spectral subspaces), and one component $\hat{\gF}_*$ that has  
fundamental group isomorphic to $\mathds{Z}$; operators in $\hat{\gF}_*$ are called strongly indefinite.   
Let now $L=\{L_t\}_{t\in [a,b]}$ be a continuous path in $\hat{\gF}$; equivalently, we can consider 
the corresponding path $Q=\{Q_t\}_{t\in [a,b]}$ of   
quadratic forms of Fredholm type on $H$. The spectral flow of $L$, or of $Q$, is the integer ${\rm sf}(L,[a,b])$, or ${\rm sf}(Q,[a,b])$,
which, intuitively, gives the net number of eigenvalues of $L$ that cross 0 in the 
positive direction as $t$ runs from $a$ to $b$. An intuitive, albeit rigorous, construction of 
${\rm sf}(L,[a,b])$ that works regardless of endpoints invertibility assumptions was carried out by Phillips in \cite{phillips}
by means of the continuous functional calculus.
According to this definition, which we shall assume, if $L$ is a path in $\hat{\gF}_+$ then
\begin{equation}\label{eq489}
{\rm sf}(L,[a,b]) = {\rm ind}\, L_a - {\rm ind}\, L_b.
\end{equation}
We shall not need the precise definition of ${\rm sf}(L,[a,b])$ in what follows, but just some basic properties that
we collect bellow and that can be found in \cite{recht}.

\medskip

\begin{enumerate}[left=4pt .. \parindent]
\item\label{invertibility} If $L_t$ is invertible for all $t$, then ${\rm sf}(L,[a,b])=0$;
\item\label{concatenationadditivity}  if $L_1 * L_2$ is the concatenation of two paths $L_1=\{(L_1)_t\}_{t\in [a,b]}$ and $L_2=\{(L_2)_t\}_{t\in[b,c]}$ that
agree for $t=b$, then
\begin{equation}\nonumber
{\rm sf}(L_1*L_2, [a,c]) = {\rm sf}(L_1,[a,b]) + {\rm sf}(L_2,[b,c]);
\end{equation}
\item\label{directsumadditivity} if $L_1$ and $L_2$ are two paths of self-adjoint Fredholm operators on Hilbert spaces $H_1$ and $H_2$, then 
\begin{equation}\nonumber
{\rm sf}(L_1\oplus L_2,[a,b]) = {\rm sf}(L_1,[a,b]) + {\rm sf}(L_2,[a,b]);
\end{equation}
\item\label{cogredience} given a path $L$ in $\hat{\gF}(H)$ and a path $M=\{M_t\}$ in $GL(H)$, then
\begin{equation}\nonumber
{\rm sf}(L,[a,b]) = {\rm sf}(M^*LM,[a,b]);
\end{equation}
\item\label{compactperturbation} 
if $L$ is closed and $K=\{K_t\}_{t\in [a,b]}$ is a
closed path of compact self-adjoint operators, then 
\begin{equation}\nonumber
{\rm sf}(L,[a,b]) = {\rm sf}(L+K,[a,b]).
\end{equation}
\end{enumerate}

Observe that, from properties \ref{invertibility}., \ref{directsumadditivity}. and \ref{cogredience}, we also have

\begin{lemma}\label{lemmakernel}
If the kernel of $L_t$ has constant dimension, then ${\rm sf}(L,[a,b])=0$.
\end{lemma}
\Proof
Let $H_1 = \ker L_a$. Since $\dim \ker L_t$ is constant, it is easy to show that we can find 
a continuous path $M=\{M_t\}$ in $GL(H)$ such that $M_t(H_1)=\ker L_t$ for all $t$. The operator 
$M_t^*L_t M_t$ then splits as a direct sum $0 \oplus L_t'$, where $0$ is the null operator on $H_1$ end 
$L_t'$ is an invertible operator on $H_1^\perp$. The result follows now from \ref{invertibility}., 
\ref{directsumadditivity}. and \ref{cogredience}. 
\qed

\subsection{Proof of Theorem \ref{maintheorem}}\label{sectionproof}

\noindent We shall only deal with the relevant case in which $Q$ is a path of strongly indefinite quadratic forms. In this case, it is clear that
$\{\overline{Q}_t\}$ are also strongly indefinite. The important fact in this regard will be the fact, used in the proof of Lemma \ref{lemmanondegenerate} below,
that the space $\hat{\mathfrak{F}}_*(V) \cap GL(V)$ is path-connected (see \cite[Lemma 3.2]{vitoriobifurcation}). 

As in \cite{vitoriobifurcation}, the proof of Theorem \ref{maintheorem} is based on the following observation, whose simple proof we shall reproduce
here for the convenience of the reader:  

\begin{lemma}\label{lemmaclosedpath}
Suppose that the path $Q$ is closed. Then ${\rm sf}(Q,[a,b]) = {\rm sf}(\overline{Q},[a,b])$.
\end{lemma}
\Proof
Let $L_t$ be the self-adjoint operator representing $Q_t$. We split $L_t$ as $L_t=N_t+K_t$, where the operators are given in block forms,
relative to the decomposition $H=V\oplus V^\perp$, as
\[
N_t = \begin{pmatrix} A_t & {\rm O} \\ {\rm O} & {\rm O} \end{pmatrix}, \quad K_t= \begin{pmatrix} {\rm O} & B_t \\ B_t^* & C_t  \end{pmatrix}.
\] 
The self-adjoint representation of $\overline{Q}_t$ is $A_t$, and since $V^\perp$ is finite dimensional, $K_t$ is compact. It thus follows from \ref{directsumadditivity}. and \ref{compactperturbation}. of Sec. \ref{sectionspectralflow} that ${\rm sf}(N,[a,b])={\rm sf}(\overline{Q},[a,b])$ 
and ${\rm sf}(N+K,[a,b]) = {\rm sf}(N,[a,b])$.    
\qed

\begin{lemma}\label{lemmanondegenerate}
If $\overline{Q}_a$ and $\overline{Q}_b$ are non-degenerate, then
\begin{equation}\nonumber
{\rm sf}(Q,[a,b]) - {\rm sf}(\overline{Q},[a,b]) = {\rm ind}\,Q_a|_{V^{\perp_{Q_a}}} 
- {\rm ind}\,Q_b|_{V^{\perp_{Q_b}}}.
\end{equation}
\end{lemma}
\Proof
The idea is to close up the path $Q=\{Q_t\}_{t\in [a,b]}$ as follows. 
Since $\overline{Q}_a$ and $\overline{Q}_b$ are non-degenerate, we have $H=V\oplus V^{\perp_{Q_a}}=V\oplus V^{\perp_{Q_b}}$. So,
for $s=a,b$, by rotating $V^{\perp_{Q_s}}$ to $V^\perp$ while keeping $V$ fixed, we can connect $Q_s$, via a path of 
quadratic forms of Fredholm type, to the quadratic form $\hat{Q}_s$ 
whose block decomposition relative to $H=V\oplus V^\perp$ is 
\begin{equation}\nonumber
\hat{Q}_s =
\begin{pmatrix}
\overline{Q}_s & {\rm O} \\
{\rm O} & Q_s^\perp
\end{pmatrix}, \quad s=a,b,
\end{equation}
for some quadratic form $Q_s^\perp$ on $V^\perp$ which is isometric to $Q_s|_{V^{\perp_{Q_s}}}$. 
Both paths have constant dimensional kernels, so by Lemma \ref{lemmakernel} their spectral flows are equal to zero.
Next we connect $\hat{Q}_b$ to $\hat{Q}_a$ by, at the same time, connecting $\overline{Q}_b$ to $\overline{Q}_a$ via a path of non-degenerate
quadratic forms on $V$ (this is possible since the space $\hat{\mathfrak{F}}_*(V)\cap GL(V)$
is path-connected) and connecting, in any way, $Q_b^\perp$ to $Q_a^\perp$ and keeping 
the null blocks unaltered. 
By additivity under direct sum, the spectral flow of this path is the sum of the spectral flows of its restrictions to 
$V$ and $V^\perp$, respectively. The first 
restriction is a path of non-degenerate quadratic forms, so has null spectral flow. The second is a path of quadratic 
forms on the finite dimensional
Hilbert space $V^\perp$ starting at $Q^\perp_b$ and ending at $Q^\perp_a$, hence its spectral flow is simply given by 
${\rm ind}\, Q^\perp_b - {\rm ind}\, Q^\perp_a = {\rm ind}\, Q_b|_{V^{\perp_{Q_b}}} - {\rm ind}\, Q_a|_{V^{\perp_{Q_a}}}$ 
(see Eq. (\ref{eq489})).  
Now, by concatenating to $Q$ the above paths we obtain a closed path $Q'=\{Q'_t\}_{t\in [a,c]}$, for some $b<c$, such that 
$Q'_t=Q_t$ for $t\in [a,b]$. From the above discussion, and from the additivity of
the spectral flow under concatenation, we have 
\[
{\rm sf}(Q',[a,c]) = {\rm sf}(Q,[a,b]) + {\rm ind}\, Q_b|_{V^{\perp_{Q_b}}} - {\rm ind}\, Q_a|_{V^{\perp_{Q_a}}}.
\]
On the other hand, since $Q'$ is a closed path, Lemma \ref{lemmaclosedpath} implies that ${\rm sf}(Q',[a,c])={\rm sf}(\overline{Q'},[a,c])$, where 
$\overline{Q'}$ is the restriction of $Q'$ to $V$. But since $\overline{Q'}_t$ is non-degenerate for $t\in[b,c]$, then
${\rm sf}(\overline{Q'},[a,c]) = {\rm sf}(\overline{Q'},[a,b])={\rm sf}(\overline{Q},[a,b])$. The result follows. \qed

\begin{lemma}\label{lemmaW}
There exists a closed finite codimensional subspace $W\subset H$ that is contained in $V$ and 
which is non-degenerate with respect to both $Q_a$ and $Q_b$.
\end{lemma}
\Proof
Let $W'\subset V$ be any closed subspace such that $V=W'\oplus \ker\overline{Q}_a$. Then $Q_a$ is non-degenerate on $W'$ and
$W'$ has finite codimension in $H$.
By the same reasoning, we can find a closed finite codimensional subspace $W''\subset H$
with $W''\subset W'$ and on which $Q_b$ is non-degenerate. Let $k$ be the codimension of $W''$ in $W'$.
From Lemma \ref{lemmadensity}, there exists a neighborhood $\cU$ of $W''$ in $G^k(W')$ such that $Q_b$ is
non-degenerate on any $W\in \cU$. Also, since $Q_a$ is non-degenerate on $W'$, the second part of Lemma \ref{lemmadensity}
applies to guarantee a $W\in \cU$ on which $Q_a$ is non-degenerate. This completes the proof.\qed

Let $W$ be as in the previous lemma. We can apply
Lemma \ref{lemmanondegenerate} twice: one time with $V$ replaced by $W$, and the other with $H$ replaced by 
$V$ (and $Q$ replaced by $Q|_V$) and $V$ replaced by $W$. We obtain, respectively, 
\begin{equation*}
\left\{
\begin{aligned}
& {\rm sf}(Q,[a,b]) - {\rm sf}(Q|_W,[a,b])  =   {\rm ind}\, Q_a|_{W^{\perp_{Q_a}}} - {\rm ind}\, Q_b|_{W^{\perp_{Q_b}}}, \\
& {\rm sf}(\overline{Q},[a,b]) - {\rm sf}(Q|_W,[a,b]) =  {\rm ind}\, Q_a|_{W^{\perp_{Q_a}}\cap V} - 
{\rm ind}\, Q_b|_{W^{\perp_{Q_b}}\cap V}.
\end{aligned}\right.
\end{equation*} 
Therefore,
\begin{equation}\label{eq4923}
{\rm sf}(Q,[a,b]) - {\rm sf}(\overline{Q},[a,b]) = \left({\rm ind}\, Q_a|_{W^{\perp_{Q_a}}} - {\rm ind}\, Q_a|_{W^{\perp_{Q_a}}\cap V}\right)
- \left( {\rm ind}\, Q_b|_{W^{\perp_{Q_b}}} - {\rm ind}\, Q_b|_{W^{\perp_{Q_b}}\cap V}\right).
\end{equation}
From now on, let $s=a,b$. Each difference of indices above can be computed by means of Eq. (\ref{eq1}) since $W^{\perp_{Q_s}}$ is finite dimensional.
Doing so, we obtain
\begin{equation}\label{eq307}
\begin{aligned}
{\rm ind}\, Q_s|_{W^{\perp_{Q_s}}} - {\rm ind}\, Q_s|_{W^{\perp_{Q_s}}\cap V}  & \;=\; 
{\rm ind}\left(Q_s\big|_{(W^{\perp_{Q_s}}\cap V)^{\perp_{Q_s}}\cap W^{\perp_{Q_s}}} \right)
+ \dim \left( (W^{\perp_{Q_s}}\cap V) \cap (W^{\perp_{Q_s}}\cap V)^{\perp_{Q_s}} \right) \\
& \quad   -\dim\left( (W^{\perp_{Q_s}}\cap V) \cap \ker Q_s|_{W^{\perp_{Q_s}}} \right).
\end{aligned}
\end{equation}
The intersections of subspaces in the right-hand side above can be simplified as follows.
On the one hand, since $H=W\oplus W^{\perp_{Q_s}}$ and $\ker Q_s\subset W^{\perp_{Q_s}}$, we have
$\ker Q_s|_{W^{\perp_{Q_s}}}=\ker Q_s$ and thus 
$(W^{\perp_{Q_s}}\cap V) \cap \ker Q_s|_{W^{\perp_{Q_s}}}=V\cap\ker Q_s$.
On the other hand, from Lemma \ref{lemmaalgebraic} we obtain  
\begin{eqnarray*}
(W^{\perp_{Q_s}}\cap V) \cap (W^{\perp_{Q_s}} \cap V)^{\perp_{Q_s}} & = &  
(W^{\perp_{Q_s}}\cap V) \cap (W\oplus V^{\perp_{Q_s}}) \\
& = & V \cap V^{\perp_{Q_s}}, \\
(W^{\perp_{Q_s}}\cap V)^{\perp_{Q_s}} \cap W^{\perp_{Q_s}} & = & (W\oplus V^{\perp_{Q_s}}) \cap W^{\perp_{Q_s}} \\
& = & V^{\perp_Q}.  
\end{eqnarray*}
The second and fourth equalities above are direct consequences of the relations 
$W\subset V$ and $W^{\perp_{Q_s}}\supset V^{\perp_{Q_s}}$ and $W\cap W^{\perp_{Q_s}}=\{0\}$.
Eq. (\ref{eq307}) can thus be rewritten as
\begin{equation}\nonumber
{\rm ind}\, Q_s|_{W^{\perp_{Q_s}}} - {\rm ind}\, Q_s|_{W^{\perp_{Q_s}}\cap V} = 
{\rm ind}\, Q_s|_{V^{\perp_{Q_s}}} + \dim(V\cap V^{\perp_{Q_s}}) - \dim (V\cap \ker Q_s).
\end{equation}
Substituting the above, for $s=a$ and $s=b$, in the right-hand side of Eq. (\ref{eq4923}) concludes the proof of Theorem \ref{maintheorem}.

\section{The Maslov index in a symplectic reduction}\label{sectionmaslov}

\noindent Throughout this section, $\mathds{H}$ shall denote a {\sl symplectic Hilbertable space}. By this we mean that $\mathds{H}$ is a separable infinite dimensional real Hilbertable space
endowed with a bounded skew-symmetric bilinear form $\omega$ which is strongly non-degenerate.

\subsection{The Fredholm Lagrangian Grassmannian and the Maslov index.}

\noindent Here we shall recall some notions regarding the Fredholm Lagrangian Grassmannian, 
$\cF \Lambda_{L_0}(\mathds{H})$, and state a definition of the Maslov index $\mu_{L_0}(\ell,[a,b])$ for continous paths $\ell:[a,b]\r \cF\Lambda_{L_0}(\mathds{H})$ in terms of
the fundamental groupoid as presented in \cite{Eidam}.

A (necessarily closed) subspace $L\subseteq \mathds{H}$ is said to be Lagrangian if $L=L^{\perp_\omega}$. The collection of all Lagrangian subspaces of $\mathds{H}$, the {\sl Lagrangian Grassmannian} of $\mathds{H}$, is denoted by $\Lambda=\Lambda (\mathds{H})$. 
An analytic Banach manifold structure for $\Lambda$ is obtained as follows. Firstly, given $L\in \Lambda$, let $\Lambda_0(L)$ be the set of all Lagrangians 
$L'$ that are complementary to $L$, i.e. $L'\cap L=\{0\}$ and $L+L'=\mathds{H}$; this is an open and dense subset of $\Lambda$. Given a pair $(L_0,L_1)$
of complementary Lagrangians, every $L\in \Lambda_0(L_1)$ is the graph of a unique bounded linear map $T:L_0\r L_1$ and the bounded bilinear form
on $L_0$ defined by $\omega[T\cdot , \cdot]|_{L_0\times L_0}$ is symmetric. This defines a chart for $\Lambda$ around $L_0$,
\begin{equation}\nonumber
\left\{
\begin{aligned}
& \varphi_{L_0,L_1}: \Lambda_0(L_1) \r B_{\rm sym}(L_0), \\
& \varphi_{L_0,L_1}(L) = \omega[T\cdot, \cdot]|_{L_0\times L_0},
\end{aligned}
\right.
\end{equation}
where $B_{\rm sym}(L_0)$ denotes the (Banach) space of all bounded symmetric bilinear forms on $L_0$. 

The {\sl Fredholm Lagrangian Grassmannian} corresponding to a given $L_0\in \Lambda$ is the subset $\cF\Lambda_{L_0} \subset \Lambda$ 
given by the $L$ such that $(L_0,L)$ is a Fredholm pair, i.e. such that $L_0\cap L$ is finite dimensional and $L_0+L$ is closed
and finite codimensional. According to \cite[Lemma 2.6]{Eidam2006}, one has

\begin{lemma}\label{lemmafredholmpair}
Given a Lagrangian $L_1$ complementary to $L_0$, a Lagrangian $L\in \Lambda_0(L_1)$ forms a 
Fredholm pair with $L_0$ if, and only if, the quadratic form $\varphi_{L_0,L_1}(L)$ is of Fredholm type.
\end{lemma}
 
It follows that $\cF \Lambda_{L_0}$ is an open subset of $\Lambda$. While the full Lagrangian Grassmannian $\Lambda$ 
is contractible (due to Kuiper's theorem),
its Fredholm counterpart $\cF\Lambda_{L_0}$ has an infinite cyclic fundamental group (see \cite{Nicolaescu}). An isomorphism 
$\pi_1(\cF\Lambda_{L_0})\cong \mathds{Z}$ is provided by the {\sl Maslov index}, whose definition can be given via the following
result from \cite{Eidam} (see also \cite{Eidam2006}).  

\begin{theorem}\label{thmmaslovfredholm}
For each $L_0\in \Lambda$, there exists a unique function assigning to each continuous path 
$\ell:[a,b]\r \cF\Lambda_{L_0}$ an integer $\mu_{L_0}(\ell,[a,b])$ that satisfies 
\begin{enumerate}[left=4pt .. \parindent]
\item $\mu_{L_0}$ is fixed-endpoint homotopy invariant;
\item $\mu_{L_0}$ is additive by concatenation; 
\item if $\ell([a,b])\subset \Lambda_0(L_1)$ for some $L_1$ complementary to $L_0$, and if $Q=\{Q_t\}$ is the path
of quadratic forms of Fredholm type on $L_0$ given by $Q_t=\varphi_{L_0,L_1}(\ell(t))$, then
\begin{equation}\nonumber
\mu_{L_0}(\ell,[a,b]) = {\rm sf}(Q,[a,b]).
\end{equation} 
\end{enumerate}
\end{theorem}

\subsection{Symplectic reduction.}\label{sectionreduction}

\noindent Let $\mathds{W}\subset \mathds{H}$ be a closed finite codimensional subspace which is {\sl coisotropic}, i.e. 
$\mathds{W}\supseteq \mathds{W}^{\perp_\omega}$. The restriction of $\omega$ to $\mathds{W}$ then descends to 
a bounded skew-symmetric and non-degenerate form $\overline{\omega}$ on
$\mathds{W}/\mathds{W}^{\perp_\omega}$ because $\ker \omega|_{\mathds{W}\times\mathds{W}} = \mathds{W}^{\perp_\omega}$. Since ${\rm codim}\,\mathds{W}<\infty$, the form 
$\omega|_{\mathds{W}\times\mathds{W}}$ is of Fredholm type. Thus $\overline{\omega}$ is also of Fredholm type and, therefore, 
it is strongly non-degenerate. The pair $(\mathds{W}/\mathds{W}^{\perp_\omega}, \overline{\omega})$ is called the symplectic reduction
of $\mathds{H}$ by $\mathds{W}$. 
 
We shall be interested in the Lagrangians $L\in \Lambda(\mathds{H})$
having null intersection with $\mathds{W}^{\perp_\omega}$. So let $\cO\subset \Lambda(\mathds{H})$ be the open set given by
\begin{equation}\nonumber
\cO=\{L \in \Lambda(\mathds{H}) \,\; | \,\; L\cap \mathds{W}^{\perp_\omega}=\{0\}\}.
\end{equation}

Let $q:\mathds{W} \r \mathds{W}/\mathds{W}^{\perp_\omega}$ be the quotient map. 

\begin{lemma}\label{lemmareduction}
We have a well-defined continuous map $\lambda:\cO \r \Lambda(\mathds{W}/\mathds{W}^{\perp_\omega})$, $\lambda(L)= q(L\cap \mathds{W})$. Moreover, 
for a fixed $L_0\in \cO$, we have $q(L\cap \mathds{W})\in \cF\Lambda_{q(L_0\cap\mathds{W})}$ if $L\in\cF\Lambda_{L_0}(\mathds{H})$.
\end{lemma}
\Proof
Given $L_0\in\cO$ we can find $L_1\in \Lambda_0(L_0)$ such that $L_1\supset \mathds{W}^{\perp_\omega}$, and hence $L_1\subset\mathds{W}$
(because $L_1^{\perp_\omega}=L_1$ and $(\mathds{W}^{\perp_\omega})^{\perp_\omega}=\mathds{W}$). 
Set $V=\mathds{W}\cap L_0$. Since $\mathds{H}=L_0\oplus L_1$ and $L_1\subset\mathds{W}$, then $\mathds{W}=V\oplus L_1$. 
It follows that $\mathds{W}/\mathds{W}^{\perp_\omega}$ is canonically $V\oplus(L_1/\mathds{W}^{\perp_\omega})$, with symplectic form
$\overline{\omega}\bigl[ (u, v+\mathds{W}^{\perp_\omega}) , (u', v'+\mathds{W}^{\perp_\omega}) \bigr]=\omega[u,v'] + \omega[v,u']$ and with 
$\lambda(L_0)$ being the summand $V$. In particular, it is clear that $\lambda(L_0)$ is Lagrangian. 
Set $\overline{L}_0=V$ and $\overline{L}_1=L_1/\mathds{W}^{\perp_\omega}$. Then $\mathds{W}/\mathds{W}^{\perp_\omega}=\overline{L}_0\oplus\overline{L}_1$
is a Lagrangian decomposition.

\smallskip

\noindent{\sl Claim.} {\sl We have $\lambda(\Lambda_0(L_1))\subset\Lambda_0(\overline{L}_1)$, and the map 
$\lambda:\Lambda_0(L_1) \r \Lambda_0(\overline{L}_1)$ corresponds, via the charts $\varphi_{L_0,L_1}$ and 
$\varphi_{\overline{L}_0,\overline{L}_1}$, to the map $B_{\rm sym}(L_0) \r B_{\rm sym}(\overline{L}_0)$,
$Q\mapsto Q|_{V=\overline{L}_0}$.}

\smallskip

\Proof 
Given $L\in\Lambda_0(L_1)$, let $T:L_0\r L_1$ be such that $L=\{u+Tu \,\; | \,\; u\in L_0\}$, so that
$\varphi_{L_0,L_1}(L)$ is the quadratic form $Q[u,v]=\omega[Tu,v]$. From $L_1\subset \mathds{W}$ we obtain
$L\cap\mathds{W}=\{u+Tu \,\; | \,\; u\in L_0\cap\mathds{W}=W\}$. It follows that $q(L\cap \mathds{W})$ is the graph
of $\overline{L}_0 \r \overline{L}_1$, $u \mapsto Tu + \mathds{W}^{\perp_\omega}$. Therefore, 
$\lambda(L)\in \Lambda_0(\overline{L}_1)$ and $\varphi_{\overline{L}_0,\overline{L}_1}(\lambda(L))$ is the 
quadratic form $\overline{Q}$ on $V=\overline{L}_0$ given by
$\overline{Q}[ u , v] = 
\overline{\omega}\bigl[ (0 , Tu + \mathds{W}^{\perp_\omega}) , (v,0) \bigr] = \omega[Tu,v] = Q[u,v].$\qed

Since $\cO$ can be covered by the open sets $\Lambda_0(L_1)$ by varying $L_1$ as above (and a fixed $L_0$), the above
claim proves the continuity of $\lambda$. The second part of the lemma follows now from Lemma \ref{lemmafredholmpair} and
the fact that $B_{\rm sym}(L_0) \r B_{\rm sym}(V)$,
$Q\mapsto Q|_V$, preserves Fredholmness of forms since $V$ has finite codimension on $L_0$. \qed

\subsection{Proof of Theorem \ref{thmreduction}.}\label{sectionproofreduction}

\noindent Since $\ell([a,b])\subset\cO$, by compactness of $[a,b]$ we can find $a=t_0<t_1<\cdots <t_k=b$ and Lagrangians $L_1$, ..., $L_k$ 
in $\Lambda_0(L_0)$ such that $\ell([t_{i-1},t_i])\subset\Lambda_0(L_i)$ and $L_i\supset\mathds{W}^{\perp_\omega}$ for all $i$. 
Fixed $i$, we shall follow the notation in the proof of Lemma \ref{lemmareduction} with $L_1$ replaced by $L_i$. So, for instance, we have
$\mathds{W}/\mathds{W}^{\perp_\omega} = \overline{L}_0\oplus \overline{L}_i$, and it was shown that $\overline{\ell}([t_{i-1},t_i])\subset \Lambda_0(\overline{L}_i)$. For each $t\in[t_{i-1},t_i]$, let 
$Q_t=\varphi_{L_0,L_i}(\ell(t)) \in B_{\rm sym}(L_0)$ and
$\overline{Q}_t = \varphi_{\overline{L}_0,\overline{L}_i}(\overline{\ell}(t))\in B_{\rm sym}(V)$. According to the proof
of Lemma \ref{lemmareduction}, $\overline{Q}_t = Q_t|_V$. Now, from 3. of Theorem \ref{thmmaslovfredholm}, we have
\begin{equation}\nonumber
\mu_{L_0}(\ell,[t_{i-1},t_i])-\mu_{\overline{L}_0}(\overline{\ell},[t_{i-1},t_i]) =
{\rm sf}(Q,[t_{i-1},t_i]) - {\rm sf}(\overline{Q},[t_{i-1},t_i]).
\end{equation} 
Since $\overline{Q}_t=Q_t|_V$, the right-hand side above is given by formula (\ref{eq2}). Therefore, by summing in $i$, the proof of
Theorem \ref{thmreduction} will be completed once we show that, for $t\in[t_{i-1},t_i]$, 
\begin{equation}\nonumber
\left\{
\begin{aligned}
& \ker Q_t = \ell(t)\cap L_0, \\
& V^{\perp_{Q_t}} = \pi\bigl( \ell(t) \cap (L_0+\mathds{W}^{\perp_\omega}) \bigr), \;\; \mbox{hence} 
\;\; V\cap V^{\perp_{Q_t}} = \pi\bigl( \ell(t)\cap(V + \mathds{W}^{\perp_\omega}) \bigr), \\
& Q_t\big|_{V^{\perp_{Q_t}}} = \mathfrak{q}_t.
\end{aligned}\right.
\end{equation}
By writing $\ell(t)$ as the graph of some $T:L_0\r L_i$, the proof of the above equalities becomes straightforward and shall
be left to the reader.

\bigskip 
\medskip


\medskip
\centerline{\headfont Data availability statement}
\medskip

No new data were created or analyzed in this study. Data sharing is not applicable to this article.

\bigskip
\centerline{\headfont References}
\medskip
\bibliography{Spectralflow_Maslovreduction}
\bibliographystyle{plain}



\end{document}